\documentclass[3p]{elsarticle}

\usepackage{amsmath}
\usepackage{amsthm}
\usepackage{amssymb}
\usepackage{bm}
\usepackage{accents}

\newtheorem{thm}{Theorem}

\newdefinition{rmk}{Remark}
\newproof{pf}{Proof}

\newcommand{\const}{\mathop{\rm const}\nolimits}

\numberwithin{equation}{section}

\journal{arXiv} 

\graphicspath{{figs/}}

\begin{document}

\begin{frontmatter}

\title{Approximate solution of the Cauchy problem for a first-order integrodifferential equation with solution derivative memory\tnoteref{label1}}
\tnotetext[label1]{The publication was financially supported by the research grant 20-01-00207 of Russian Foundation for Basic Research}

\author{P.N. Vabishchevich\corref{cor1}\fnref{lab1,lab2}}
\ead{vabishchevich@gmail.com}
\cortext[cor1]{Correspondibg author.}

\address[lab1]{Nuclear Safety Institute, Russian Academy of Sciences,
              52, B. Tulskaya, 115191 Moscow, Russia}

\address[lab2]{North-Caucasus Center for Mathematical Research, North-Caucasus Federal University, 
                   1, Pushkin Street, 355017 Stavropol, Russia}

\begin{abstract}
We consider the Cauchy problem for a first-order evolution equation with memory in a finite-dimensional Hilbert space when the integral term is related to the time derivative of the solution.
The main problems of the approximate solution of such nonlocal problems are due to the necessity to work with the approximate solution for all previous time moments.
We propose a transformation of the first-order integrodifferential equation to a system of local evolutionary equations.
We use the approach known in the theory of Voltaire integral equations with an approximation of the difference kernel by the sum of exponents. 
We formulate a local problem for a weakly coupled system of equations with additional ordinary differential equations.
We have given estimates of the stability of the solution by initial data and the right-hand side for the solution of the corresponding Cauchy problem.
The primary attention is paid to constructing and investigating the stability of two-level difference schemes, which are convenient for computational implementation.
The numerical solution of a two-dimensional model problem for the evolution equation of the first order, when the Laplace operator conditions the dependence on spatial variables, is presented.
\end{abstract}


\end{frontmatter}

\section{Introduction}\label{sec:1}

In applied mathematical modeling of nonstationary processes, parabolic and hyperbolic equations are most widely used \cite{LionsBook,evans2010partial}, the boundary value problems for which are most well studied.  
Recently, more and more attention has been paid to equations that partially inherit both the properties of parabolic and hyperbolic equations. An example is the evolutionary integrodifferential equations \cite{GripenbergBook1990,pruss2013evolutionary}.
The most important feature of such equations is their nonlocality, the dependence of the solution on the entire prehistory of the process. 

We can distinguish two classes of evolutionary integrodifferential equations with memory.
The first one is characterized by the nonlocality of the solution when the integrand includes the solution itself. 
The second class of equations is characterized by the nonlocality of the time derivative of the solution. 
Such mathematical models are typical when considering dynamic viscoelastic processes \cite{christensen1982theory,marques2012computational}.
Nowadays, the integrodifferential equations with time derivative solution memory are often associated with time-fractional equations \cite{kochubei2019equations}.  
The noted division of problems for equations with memory is rather conventional. For problems with a difference kernel, we can pass from one type of nonlocality to another when introducing another difference kernel.
Because of this, the same viscoelastic equations can (see, e.g., \cite{LionsBook}) be written as integrodifferential equations with solution memory rather than solution time derivative.

The approximate solution of boundary value problems for equations with memory is carried out using standard finite-element or finite-volume approximations on the space \cite{KnabnerAngermann2003,QuarteroniValli1994}.
We arrive at the Cauchy problem for operator equations with memory in the corresponding finite-dimensional Hilbert space.
Our primary focus should be on the problems of choosing approximations over time.
For problems with solution memory, it is natural to focus \cite{ChenBook1998} on the use of quadratures for the integral term and the usual two-level approximations of the time derivative. 
Such studies of the implicit Euler scheme and the Crank-Nicholson scheme are done, e.g., in \cite{mclean1993numerical,mclean1996discretization}. 
When considering integrodifferential equations with time derivative solution memory, the corresponding quadrature formulas are used for the integral term with time derivative. 
Various variants for time-fractional equations are actively discussed in the literature \cite{baleanu2012fractional}. 

Standard computational algorithms for approximate solutions of the Cauchy problem for integrodifferential equations with memory involve the need to work with the solution at the last moments. 
Some possibilities for reducing computational work for problems with fractional time derivatives are discussed in \cite{diethelm2020good}.
For us, the approaches with the transition from a nonlocal to a local problem, when the memory requirements, in particular, are significantly reduced, will be of most interest.

For the Volterra integral equations, the well-known (see, e.g., \cite{linz1985analytical}) transition to more computationally simple problems is provided by choosing particular approximations of the difference kernel.
Let us single out as the most promising approximation of the difference kernel by the sum of exponents. 
For equations with fractional time derivatives, this approach has been used in various variants (see, e.g., \cite{jiang2017fast,zhang2021exponential}) since \cite{baffet2017kernel}.  

When approximating the kernel by a sum of exponents, we arrive at a system of local weakly coupled evolution equations.
We considered the possibilities of such an approach for approximate solution of the Cauchy problem for an integrodifferential first-order equation with solution memory in \cite{vabMemory}. 
A similar study for problems with time derivative memory is carried out in the present paper.

We give brief information about the content of the paper. Section 2 formulates a Cauchy problem for the integrodifferential Volterra equation with a positive definite self-adjoint operator in a real finite-dimensional Hilbert space when the time derivative is inherited.
The difference kernel is approximated by the sum of exponents in Section 3. 
The transformation of the nonlocal problem for the equation with memory into a local system of equations is performed.  
Appropriate a priori estimates for the solution of the Cauchy problem are obtained.
In Section 4, two-level difference schemes with a convenient computational realization have been proposed and investigated for their stability.
The results of the numerical solution of the two-dimensional model problem are presented in Section 5.
Conclusions on the conducted research are given in Section 7.

\section{Problem formulation}\label{sec:2}

We consider a Cauchy problem for an evolutionary equation with derivative solution memory in a real finite-dimensional Hilbert space $H$.
The function $u(t)$ satisfies an integrodifferential equation with a difference kernel
\begin{equation}\label{2.1}
 B\frac{d u}{d t} + \int_{0}^{t} k(t-s) C \frac{d u}{d s}(s) d s + A u = f(t),
 \quad t > 0 ,  
\end{equation} 
and the initial condition
\begin{equation}\label{2.2}
 u(0) = u_0 .
\end{equation} 
Linear constants (independent of $t$) of operators $A, B, C$ are self-adjoint and positive definite:
\begin{equation}\label{2.3}
  A = A^* \geq \nu_A  I,   \quad \nu_A   > 0 ,   
  \quad B = B^* \geq \nu_B I,   \quad \nu_B   > 0 ,   
  \quad C = C^* \geq \nu_C  I,   \quad \nu_C   > 0 , 
\end{equation} 
where $I$ is the identity operator in $H$.
We will use the usual notations $(\cdot, \cdot)$ and $\|\cdot\|$ for the scalar product and norm in $H$.
For a self-adjoint and positive operator $D$, the Hilbert space $H_D$ is defined
with scalar product and norm $(u, v)_D = (D u, v), \ \| u \|_D = (u, v)_D ^ {1/2} $.

As when considering equations with solution memory \cite{mclean1993numerical,mclean1996discretization}, the kernel $k(t)$ is assumed to be real-valued and positive definite (convolution kernels of positive type \cite{gripenberg1990volterra}).
In this case, for each $T > 0$ the kernel $k(t)$ belongs to $L_1(0, T)$ and satisfies the inequality 
\begin{equation}\label{2.4}
 \int_{0}^{T} \psi(t) \int_{0}^{t} k(t-s) \psi(s) d s \, d t \geq 0,
 \quad \psi \in C[0,T] .
\end{equation} 
Note also \cite{halanay1965asymptotic} a sufficient condition of positive definite kernel $k(t)$:
\begin{equation}\label{2.5}
 k(t) \geq 0,
 \quad \frac{d k}{d t}  (t) \leq 0,
 \quad \frac{d^2 k}{d t^2}  (t) \geq 0,
 \quad t > 0 . 
\end{equation} 

In our study, we focus on obtaining a system of local evolution equations, the Cauchy problem for which gives an approximate solution to the problem (\ref{2.1}), (\ref{2.2}),
When investigating the stability of difference approximations in time, the following statement is our guideline.

\begin{thm}\label{t-1}
Let the operators $A, B, C$ satisfy conditions (\ref{2.3}) and let $k(t)$ be a positive definite kernel. 
Then, for the solution of the problem (\ref{2.1}), (\ref{2.2}), the stability estimate for the initial data and the right-hand side 
\begin{equation}\label{2.6}
 \|u(t)\|^2_A \leq \|u_0\|^2_A + \frac{1}{2} \int_{0}^{t} \|f(s)\|^2_{B^{-1}} d s,
 \quad t > 0 , 
\end{equation} 
is valid.
\end{thm}

\begin{pf}
Let us multiply equation (\ref{2.1}) scalarly in $H$ by $du(t)/dt$ and, given the positive definiteness of operators $A, B, C$, we obtain
\[
 \frac{1}{2} \frac{d}{d t} \|u(t)\|_A^2 + \int_{0}^{t} k(t-s)  \Big( C^{1/2} \frac{d u}{d s}(s),  C^{1/2} \frac{d u}{d t}(t) \Big) d s 
 \leq  \frac{1}{4} \|f (t)\|^2_{B^{-1}} .
\] 
By integration over $(0, T)$, this yields 
\[
 \frac{1}{2} \big(\|u(T)\|_A^2 - \|u(0)\|_A^2 \big) 
 +  \int_{0}^{T} \int_{0}^{t} k(t-s)  \Big(C^{1/2} \frac{d u}{d s}(s), C^{1/2}\frac{d u}{d t}(t) \Big)  d s \ d t 
 \leq  \frac{1}{4} \int_{0}^{T}  \|f(t)\|^2_{B^{-1}} d t.
\]  
Given (\ref{2.4}) and the initial condition (\ref{2.2}), we have a provable estimate (\ref{2.6}). 
\end{pf}

The nonlocal term in (\ref{2.1}) becomes local in two important cases: when the kernel $k(t)$ is constant and when the kernel is a $\delta-$function. 
We can distinguish such terms separately:
\begin{equation}\label{2.7}
 k(t) \rightarrow  \gamma_1  + \gamma_2 \delta(t) + k(t), 
 \quad \gamma_1 >  0,  
 \quad \gamma_2 >  0 .  
\end{equation} 
This corresponds to the transition 
\begin{equation}\label{2.8}
 B \rightarrow B + \gamma_1 C,
 \quad A \rightarrow A + \gamma_2 C,
 \quad f(t) \rightarrow f(t) + \gamma_2 C u_0 ,
\end{equation} 
in equation (\ref{2.1}).
Thus we remain in the class of problems under consideration (\ref{2.1})--(\ref{2.3}).

\section{Transformation to a local problem}\label{sec:3}

Computationally, the approaches with constructing an approximate solution of the nonlocal problem (\ref{2.1}), (\ref{2.2}) with the memory of the solution derivative in time using the solutions of some local evolutionary problems are of the most significant interest. 
They can be created based on the approximation of the kernel by a sum of exponents.

The kernel $k(t)$ is approximated by the function $\accentset{\sim }{k}(t)$, which has the form
\begin{equation}\label{3.1}
 \accentset{\sim }{k}(t) = \sum_{i=1}^{m} a_i \exp(-b_i t),
 \quad t \geq 0 . 
\end{equation} 
For the coefficients $a_i, b_i, \ i = 1,2,\ldots, m,$ assumptions are made 
\begin{equation}\label{3.2}
 a_i > 0,
 \quad b_i > 0,
 \quad   i = 1,2,\ldots, m .
\end{equation} 
Given the conditions (\ref{2.5}) under the constraints (\ref{3.2}), the kernel $\accentset{\sim }{k}(t)$ is a positive definite kernel. 

We denote by $v(t)$ the approximate solution of the problem (\ref{2.1}), (\ref{2.2}).
It is defined as the solution of the Cauchy problem
\begin{equation}\label{3.3}
 B \frac{d v}{d t} + \int_{0}^{t} \accentset{\sim }{k}(t-s)  C \frac{d v}{d s}(s) d s  + A v = f (t),
 \quad t > 0 ,  
\end{equation} 
\begin{equation}\label{3.4}
 v(0) = u_0 .
\end{equation} 

To pass from the nonlocal problem (\ref{3.3}), (\ref{3.4}) to the local one, we introduce \cite{linz1985analytical,vabMemory} functions
\[
 v_i(t) = \int_{0}^{t} \exp(-b_i (t-s)) \frac{d v}{d s}(s) d s ,
 \quad i = 1,2, \ldots, m . 
\] 
Given the approximation (\ref{3.1}), the equation (\ref{3.3}) is written as
\begin{equation}\label{3.5}
 B \frac{d v}{d t} + \sum_{i=1}^{m} a_i C v_i + A v = f (t) .
\end{equation} 
For the auxiliary functions $v_i(t), \ i = 1,2, \ldots, m,$ we have the equations
\begin{equation}\label{3.6}
 \frac{d v_i}{d t} + b_i v_i - \frac{d v}{d t} = 0 ,
 \quad i = 1,2, \ldots, m . 
\end{equation} 
The system of the equations (\ref{3.5}), (\ref{3.6}) is supplemented by the initial conditions
\begin{equation}\label{3.7}
 v(0) = u_0,
 \quad v_i(0) = 0,
 \quad i = 1,2, \ldots, m .  
\end{equation} 

It is convenient to rewrite the equation (\ref{3.5}) in a slightly different form.
Substitution 
\[
  v_i = \frac{1}{b_i} \frac{d v}{d t} - \frac{1}{b_i} \frac{d v_i}{d t} ,
 \quad i = 1,2, \ldots, m ,  
\] 
in the equation (\ref{3.5}) gives
\begin{equation}\label{3.8}
 \Big (B + \sum_{i=1}^{m} \frac{a_i}{b_i} C \Big ) \frac{d v}{d t}  -  \sum_{i=1}^{m} \frac{a_i}{b_i} C \frac{d v_i}{d t} + A v = f (t) .
\end{equation}

\begin{thm}\label{t-2}
Let the operators $A, B, C$ satisfy the conditions (\ref{2.3}). Then
for the solution of the problem (\ref{3.2}), (\ref{3.6})--(\ref{3.8}) the stability estimation on the initial data and the right-hand side 
\begin{equation}\label{3.9}
 \|v(t)\|^2_A + \sum_{i=1}^{m} a_i \|v_i(t)\|^2_{C} \leq \|u_0\|^2_A + \frac{1}{2}  \int_{0}^{t} \|f(s)\|^2_{B^{-1}} d s ,
 \quad t > 0 , 
\end{equation} 
is valid.
\end{thm}

\begin{pf}
We multiply the equation (\ref{3.8}) by $d v(t) / dt$, and the separate equation (\ref{3.6}) for $i = 1,2, \ldots, m$ by  $a_i b_i^{-1} C d v_i(t)/d t$. This gives equalities
\[
 \Big \|\frac{d v}{d t} \Big \|^2_B + 
 \frac{1}{2}  \frac{d}{d t} \|v\|^2_A + \sum_{i=1}^{m} \frac{a_i}{b_i} \Big (C \frac{d v}{d t}, \frac{d v}{d t} \Big) 
 - \sum_{i=1}^{m} \frac{a_i}{b_i} \Big (C \frac{d v_i}{d t}, \frac{d v}{d t} \Big) = \Big (f, \frac{d v}{d t} \Big ),
\] 
\[
 \frac{a_i}{b_i} \Big (C \frac{d v_i}{d t}, \frac{d v_i}{d t} \Big) - \frac{a_i}{b_i} \Big (C \frac{d v}{d t}, \frac{d v_i}{d t} \Big) 
 + a_i \frac{1}{2} \frac{d}{d t} \|v_i\|_{C}^2 = 0 ,
 \quad i = 1,2, \ldots, m .  
\] 
Adding them, we get
\[
 \Big \|\frac{d v}{d t} \Big \|^2_B + \frac{1}{2}  \frac{d}{d t} \Big ( \|v\|^2_A + \sum_{i=1}^{m} a_i  \|v_i\|_{C}^2 \Big ) 
 + \sum_{i=1}^{m} \frac{a_i}{b_i} \Big \|\frac{d v}{d t} - \frac{d v_i}{d t} \Big \|^2_C = \Big (f, \frac{d v}{d t} \Big ) .
\] 
This gives the inequality
\[
 \frac{d}{d t} \Big ( \|v\|^2_A + \sum_{i=1}^{m} a_i  \|v_i\|_{C}^2 \Big )  
 \leq \frac{1}{2}  \|f(t)\|^2_{B^{-1}} .
\] 
From this follows the provable estimate (\ref{3.9}).
\end{pf}

It is convenient for us to write the system (\ref{3.6}), (\ref{3.8}) as one first-order equation for vector quantities.
Define the vector 
$\bm v = \{v, v_1, \ldots, v_m \} $ and $\bm f  = \{f, 0, \ldots, 0 \} $, 
and from (\ref{3.6}), (\ref{3.8}), we get to the Cauchy problem
\begin{equation}\label{3.10}
 \bm B \frac{d \bm v}{d t} + \bm A \bm v = \bm f  ,
\end{equation} 
\begin{equation}\label{3.11}
 \bm v(0) = \bm v_0 ,
\end{equation} 
where $\bm v_0 = \{u_0, 0, \ldots, 0 \}$.
For the operator matrices $\bm B$ and $\bm A$,  we have the representation
\begin{equation}\label{3.12}
\bm B = \left (\begin{array}{cccc}
  B + {\displaystyle \sum_{i=1}^{m} \frac{a_i}{b_i} C  } & - {\displaystyle \frac{a_1}{b_1} C}  & \cdots &  - {\displaystyle \frac{a_m}{b_m} C}   \vspace{2mm} \\
  - {\displaystyle \frac{a_1}{b_1} C }  &  {\displaystyle \frac{a_1}{b_1} C } & \cdots &  0 \\
  \cdots  & \cdots & \cdots &  0 \\
  - {\displaystyle \frac{a_m}{b_m} C } &  0 & \cdots &  {\displaystyle \frac{a_m}{b_m} C } \\
\end{array}
 \right ) ,
 \quad \bm A = \mathrm{diag} \, \big (A, a_1 C, \ldots, a_m C \big) .
 \quad 
\end{equation} 

The problem (\ref{3.10}), (\ref{3.11}) we consider on the direct sum of spaces $\bm H = H \oplus \ldots  \oplus H$,
when for $\bm v, \bm w \in \bm H$,  the scalar product and norm are determined by the expressions
\[
 (\bm v, \bm w) = (v, w) + \sum_{1=1}^{m} (v_i, w_i),
 \quad \|\bm v\| =  (\bm v, \bm v)^{1/2} .
\]
Given the conditions (\ref{2.3}) and (\ref{3.2}), we obtain
\begin{equation}\label{3.13}
 \bm B = \bm B^* \geq 0,
 \quad \bm A = \bm A^* > 0 .
\end{equation} 

To prove the estimate (\ref{3.9}), we multiply scalarly in $\bm H$ the equation (\ref{3.10}) $d \bm v / dt$.
Given the properties (\ref{3.13}), this gives
\[
 \Big (\bm B \frac{d \bm v}{d t}, \frac{d \bm v}{d t} \Big ) + 
 \frac{1}{2} \frac{d}{d t} \|\bm v\|^2_{\bm A} = \Big (\bm f , \frac{d \bm v}{d t} \Big ) .  
\]
Considering 
\[
 \Big (\bm B \frac{d \bm v}{d t}, \frac{d \bm v}{d t} \Big ) \geq \Big (B \frac{d v}{d t}, \frac{d v}{d t} \Big ) ,
 \quad \Big (\bm f , \frac{d \bm v}{d t} \Big ) \leq \Big (B \frac{d v}{d t}, \frac{d v}{d t} \Big ) + \frac{1}{4} (B^{-1} f, f) ,
\] 
we have 
\begin{equation}\label{3.14}
 \|\bm v(t)\|^2_{\bm A} \leq \|\bm v_0\|^2_{\bm A} +  \frac{1}{2} \int_{0}^{t} \|f(s)\|_{B^{-1}} d s .
\end{equation} 
In our case
\[
 \|\bm v(t)\|^2_{\bm A} = \|v(t)\|^2_A + \sum_{i=1}^{m} a_i \|v_i(t)\|_{C}^2 , 
 \quad \|\bm v_0\||^2_{\bm A} = \|u_0\|^2_A ,
\] 
so the inequality (\ref{3.14}) gives the estimate (\ref{3.9}).

Instead of approximating (\ref{3.1}), we can investigate the slightly more general case where similarly (\ref{2.7}) 
\begin{equation}\label{3.15}
 \accentset{\sim }{k}(t) = \gamma_1 + \gamma_2 \delta(t) + \sum_{i=1}^{m} a_i \exp(-b_i t),
 \quad t \geq 0 . 
\end{equation} 
The transition to the considered case is provided by (\ref{2.8}).  

\section{Two-level difference schemes}\label{sec:4} 

In the approximate solution of the Cauchy problem (\ref{3.10}), (\ref{3.11}),  implicit time approximations are often used.
In this case, we have unconditionally stable schemes. 
We will use, for simplicity, a uniform grid in time with step $\tau$ and let  $y^n=y(t^n), \ t^n = n\tau$, $n =0, 1, \ldots$. 
We consider a two-level scheme with the weight $\sigma = \const \in (0,1]$, when 
\begin{equation}\label{4.1}
 \bm B \frac{\bm y^{n+1} - \bm y^{n}}{\tau } + \bm A \bm y^{n+\sigma} = \bm f^{n+\sigma},
 \quad n = 0,1,\ldots,
\end{equation} 
\begin{equation}\label{4.2}
 \bm y^0 = \bm v_0 ,
\end{equation} 
when using the notation
\[
 \bm y^{n+\sigma} = \sigma \bm y^{n+1} + (1-\sigma ) \bm y^{n} ,
 \quad   \bm y^{n} = \{y^{n}, y_1^{n}, \ldots, y_m^{n}\} .
\] 
For the right-hand side and the initial condition, we have
\[
 \bm f^{n+\sigma} = \{f^{n+\sigma}, 0, \ldots, 0 \},
 \quad   \bm v^{0} = \{u_0, 0, \ldots, 0 \} . 
\] 

The difference scheme (\ref{4.1}), (\ref{4.2}) approximates the problem (\ref{3.10}), (\ref{3.11}) with sufficient smoothness of the solution $\bm v(t)$ with the first order in $\tau$ for $\sigma \neq 0.5$ and with the second order for $\sigma = 0.5$ (Crank-Nicolson scheme).
To study the stability of two-level schemes, we can use the results of the theory of stability (correctness) of operator-difference schemes \cite{SamarskiiTheory,SamarskiiMatusVabischevich2002}.

\begin{thm}\label{t-3}
The two-level scheme (\ref{2.3}), (\ref{3.12}), (\ref{4.1}), (\ref{4.2}) 
is unconditionally stable for $\sigma \geq 0.5$.
Under these constraints, for an approximate solution to the problem (\ref{3.10}), (\ref{3.11}), the a priori estimate 
\begin{equation}\label{4.3}
 \|\bm y^{n+1}\|^2_{\bm A} \leq \|u_0\|^2_{A}  + \frac{1}{2} \sum_{k=0}^{n}\tau  \|f^{k+\sigma}\|^2_{B^{-1}} ,
 \quad n = 0,1,\ldots ,
\end{equation} 
holds.
\end{thm}

\begin{pf}
Let us write the equation (\ref{4.1}) as
\[
 \Big (\bm B + \Big (\sigma - \frac{1}{2} \Big ) \tau  \bm A \Big ) \frac{\bm y^{n+1} - \bm y^{n}}{\tau } 
 + \bm A \frac{\bm y^{n+1} + \bm y^{n}}{2} = \bm f^{n+\sigma} .
\] 
By multiplying this equation by $2 (\bm y^{n+1}- \bm y^{n})$, given $\sigma \geq 0.5$, we obtain
\begin{equation}\label{4.4}
 2 \tau \Big (\bm B \frac{\bm y^{n+1} - \bm y^{n}}{\tau }, \frac{\bm y^{n+1} - \bm y^{n}}{\tau } \Big ) 
 + \|\bm y^{n+1}\|^2_{\bm A} - \|\bm y^{n}\|^2_{\bm A} \leq 2 \tau  \Big (\bm f^{n+\sigma}, \frac{\bm y^{n+1} - \bm y^{n}}{\tau } \Big ).
\end{equation}  
Since
\[
 \Big (\bm B \frac{\bm y^{n+1} - \bm y^{n}}{\tau }, \frac{\bm y^{n+1} - \bm y^{n}}{\tau } \Big )
 \geq \Big (B \frac{y^{n+1} - y^{n}}{\tau }, \frac{y^{n+1} - y^{n}}{\tau } \Big ) ,
\] 
\[
 \Big (\bm f^{n+\sigma}, \frac{\bm y^{n+1} - \bm y^{n}}{\tau } \Big ) 
 = \Big (f^{n+\sigma}, \frac{y^{n+1} - y^{n}}{\tau } \Big ) 
 \leq \Big (B \frac{y^{n+1} - y^{n}}{\tau }, \frac{y^{n+1} - y^{n}}{\tau } \Big ) 
 + \frac{1}{4} \big (B^{-1}f^{k+\sigma}, f^{k+\sigma} \big ) ,
\] 
then from the inequality (\ref{4.4}) follows the estimate
\[
 \|\bm y^{n+1}\|^2_{\bm A} - \|\bm y^{n}\|^2_{\bm A} \leq \frac{1}{2} \tau \|f^{k+\sigma}\|^2_{B^{-1}} .
\] 
This leads us to the estimate (\ref{4.3}), which acts as the grid analogue of the estimate (\ref{3.14}).
\end{pf}

We will write a scheme with weights (\ref{4.1}), (\ref{4.2}) for the individual components.
It corresponds to the case where a difference scheme 
\begin{equation}\label{4.5}
 \Big (B + \sum_{i=1}^{m} \frac{a_i}{b_i} C \Big ) \frac{y^{n+1} - y^{n}}{\tau } 
 -  \sum_{i=1}^{m} \frac{a_i}{b_i} C \frac{y_i^{n+1} - y_i^{n}}{\tau}  + A y^{n+\sigma} = f ^{n+\sigma },
\end{equation} 
\begin{equation}\label{4.6}
 \frac{y_i^{n+1} - y_i^{n}}{\tau} + b_i y_i^{n+\sigma } - \frac{y^{n+1} - y^{n}}{\tau } = 0 ,
 \quad i = 1,2, \ldots, m ,
 \quad n = 0,1, \ldots ,  
\end{equation} 
\begin{equation}\label{4.7}
 y^{0} = u_0,
 \quad y_i^{0} = 0,
 \quad i = 1,2, \ldots, m , 
\end{equation} 
is used for the approximate solution of the problem (\ref{3.6})-(\ref{3.8}).
The inequality (\ref{4.3}) implies the a priori estimate
\begin{equation}\label{4.8}
 \|y^{n+1}\|^2_A + \sum_{i=1}^{m} a_i \|y_i^{n+1}\|_{C}^2  \leq \|u_0\|^2 + \frac{1}{2} \sum_{k=0}^{n}\tau \|f^{n+\sigma }\|^2_{B^{-1}} ,
 \quad n = 0,1, \ldots ,
\end{equation} 
for an approximate solution of the problem (\ref{4.5})--(\ref{4.7}).
The estimate (\ref{4.8}) is a difference analogue of the estimate (\ref{3.9}) for solving the differential problem (\ref{3.6})--(\ref{3.8}).

The problem of computational realization deserves special attention when solving nonlocal problems.
In the case of the equation (\ref{4.6}), we have
\begin{equation}\label{4.9}
 y_i^{n+1} = \frac{1}{1 + \sigma b_i \tau } y^{n+1} + \chi_i^{n} ,
 \quad \chi_i^{n} =  \frac{1 }{1 + \sigma b_i \tau } \Big ( \big (1 - (1-\sigma) b_i \tau \big ) y_i^{n}  - y^{n} \Big ) ,
 \quad i = 1,2,\ldots, m .
\end{equation}  
Substituting this into the equation (\ref{4.5}) gives the equation 
\begin{equation}\label{4.10}
 \big (B + \sigma \tau (\mu C +  A) \big ) y^{n+1} = \chi^{n} 
\end{equation} 
for finding $y^{n+1}$.
For the coefficient $\mu$ and the right-hand side, we have
\[
 \mu = \sum_{i=1}^{m} \frac{a_i }{1 + \sigma b_i \tau } ,
 \quad \chi^{n} = \tau f^{n+\sigma } + (B - (1-\sigma)\tau A) y^{n} +  \sum_{i=1}^{m} \frac{a_i}{b_i} C ( y^{n} -  y_i^{n} + \chi_i^{n}) .
\] 
Thus, the transition to a new $n+1$ level in time is provided by solving the standard problem (\ref{4.10}) for $y^{n+1}$ and calculating the auxiliary quantities $y_i^{n+1}, \ i = 1,2,\ldots, m,$ according to (\ref{4.9}).
The computational complexity of the approximate solution of the nonlocal problem under consideration (\ref{2.1}), (\ref{2.5}), (\ref{2.5}) is not much greater than that of the local problem. It is necessary to operate additionally with the solutions of $m$ simple auxiliary local evolution problems in explicit calculations of their solutions at a new time level.

\section{Numerical experiments }\label{sec:5} 

We will illustrate the possibilities of the proposed computational algorithms by the results of the numerical solution of a model two-dimensional problem.
We will assume that the computational domain is a unit square:
\[
 \Omega = \{ \bm x  \ | \ \bm x = (x_1,x_2), \ 0 < x_d  < 1 , \ d  = 1,2 \} ,
\]
with boundary $\partial \Omega$.
The function $w(\bm x,t)$ satisfies the equation
\[
 \frac{\partial w}{\partial t} + c \int_{0}^{t} k(t-s) \frac{\partial w}{\partial s} (\bm x,s) d s 
 - \triangle w = 0,
 \quad \bm x \in \Omega ,
 \quad 0 < t \leq T , 
\] 
with a nonnegative numerical parameter $c > 0$. 
The boundary and initial conditions have the form
\[
  w(\bm x,t) = 0,
  \quad \bm x \in \partial \Omega ,
 \quad 0 < t \leq T ,  
\] 
\[
 w (\bm x,0) = u_0(\bm x) ,
 \quad \bm x \in \Omega .
\] 

For the numerical solution of this boundary value problem, we will use standard difference approximations in space \cite{SamarskiiTheory}.
We will introduce in the region $\Omega$ a uniform rectangular grid
\[
\overline{\omega}  = \{ \bm{x} \ | \ \bm{x} =\left(x_1, x_2\right), \quad x_d =
i_d  h_d , \quad i_d  = 0,1,...,N_d ,
\quad N_d  h_d  = 1 , \quad  d  = 1,2 \} ,
\]
where $\overline{\omega} = \omega \cup \partial \omega$, 
$\omega$ is the set of internal mesh nodes, and $\partial \omega$ is the set of boundary mesh nodes.
For grid functions $w(\bm x)$ such that $w(\bm x) = 0, \ \bm x \notin \omega$, we define the Hilbert space
$H = L_2 (\omega)$, in which the scalar product and norm are
\[
(w, u) = \sum_{\bm x \in  \omega} w(\bm{x}) u(\bm{x}) h_1 h_2,  \quad 
\| w \| =  (w,w)^{1/2}.
\]
For $w(\bm x) = 0, \ \bm x \notin \omega$, we define the grid Laplace operator $- A$ on the usual five-point stencil:
\[
  \begin{split}
  A w = & -
  \frac{1}{h_1^2} (w(x_1+h_1,x_2) - 2 w(\bm{x}) + w(x_1-h_1,x_2)) \\ 
  & - \frac{1}{h_2^2} (w(x_1,x_2+h_2) - 2 w(\bm{x}) + w(x_1,x_2-h_2)), 
  \quad \bm{x} \in \omega . 
 \end{split} 
\] 
On sufficiently smooth functions, the operator $A$ approximates the differential operator
$-\triangle$ with an error $\mathcal{O} \left(|h|^2\right)$, $|h|^2 = h_1^2+h_2^2$. 
This grid operator (see, e.g., \cite{SamarskiiTheory}) is self-adjoint and positive definite in $H$.
Approximation in space leads us to the problem (\ref{2.1}), (\ref{2.2}) in which $B = I, \ C = c I$.  

In the numerical results presented below, the kernel is
\begin{equation}\label{5.1}
 k(t) = \frac{1}{\Gamma (1-\alpha)} t^{-\alpha} \exp\big(- \delta t \big),
 \quad 0 < \alpha   < 1 ,
 \quad \delta \geq 0 ,
\end{equation} 
where $\Gamma(\cdot)$ is the Gamma function.
For $\delta = 0$, such a kernel is associated with the Caputo fractional derivative.

Constructing approximations (\ref{3.1}) is an independent task.
In the approximation of nonlinear functions, the most widespread \cite{braess1986nonlinear} are rational approximations.
Theoretical and practical developments in rational approximation can be used when approximating the difference kernel by the sum of exponentials.

Define the function $K(s)$ for real $s \geq 0$ as the Laplace transform of the kernel $k(t)$:
\[
 K(s) = \int_{0}^{\infty } k(t) \exp(-st) d t . 
\] 
Let's assume that there is a rational approximation for $K(s)$ in the form
\begin{equation}\label{5.2}
 \accentset{\sim }{K}(s) = \sum_{i=1}^{m} \frac{a_i}{b_i + s} 
\end{equation} 
with coefficients $a_i, b_i, \ i = 1,2,\ldots, m,$ satisfying (\ref{3.2}).
For the original $\accentset{\sim }{k}(t)$ we obtain a representation in the form (\ref{3.1}).
Thus, the rational approximation (\ref{5.2}) of the Laplace transform $K(s)$ gives the approximation (\ref{3.1}) of the kernel by the sum of exponents.
With the case (\ref{3.15}) we associate the rational approximation
\begin{equation}\label{5.3}
 \accentset{\sim }{K}(s) = \gamma_1 \frac{1}{s} + \gamma_2 +  \sum_{i=1}^{m} \frac{a_i}{b_i + s} .
\end{equation} 
For the kernel (\ref{5.1}), we have
\[
 K(s) = (s + \delta)^{\alpha - 1} . 
\] 
The best uniform rational approximations to real scalar functions in the setting of zero defect \cite{hofreither2021algorithm} and the software developed by Clemens Hofreither. The \textsf{baryrat} open-source Python package (\textsf{https://github.com/c-f-h/baryrat}) are used to construct approximations (\ref{5.3}) with $\gamma_1 \equiv 0$. 

\begin{figure}
\centering
\includegraphics[width=0.49\linewidth]{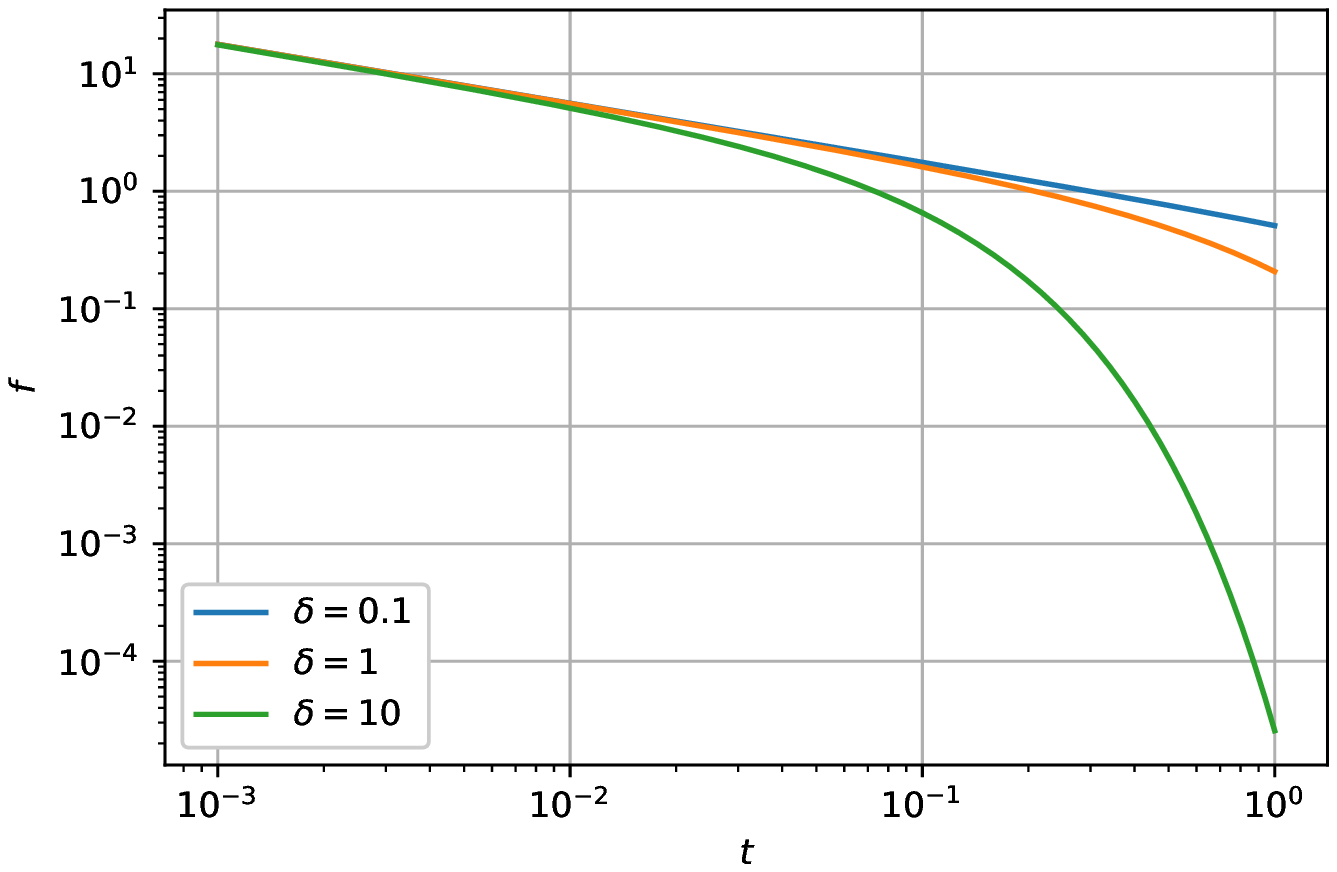} 
\includegraphics[width=0.49\linewidth]{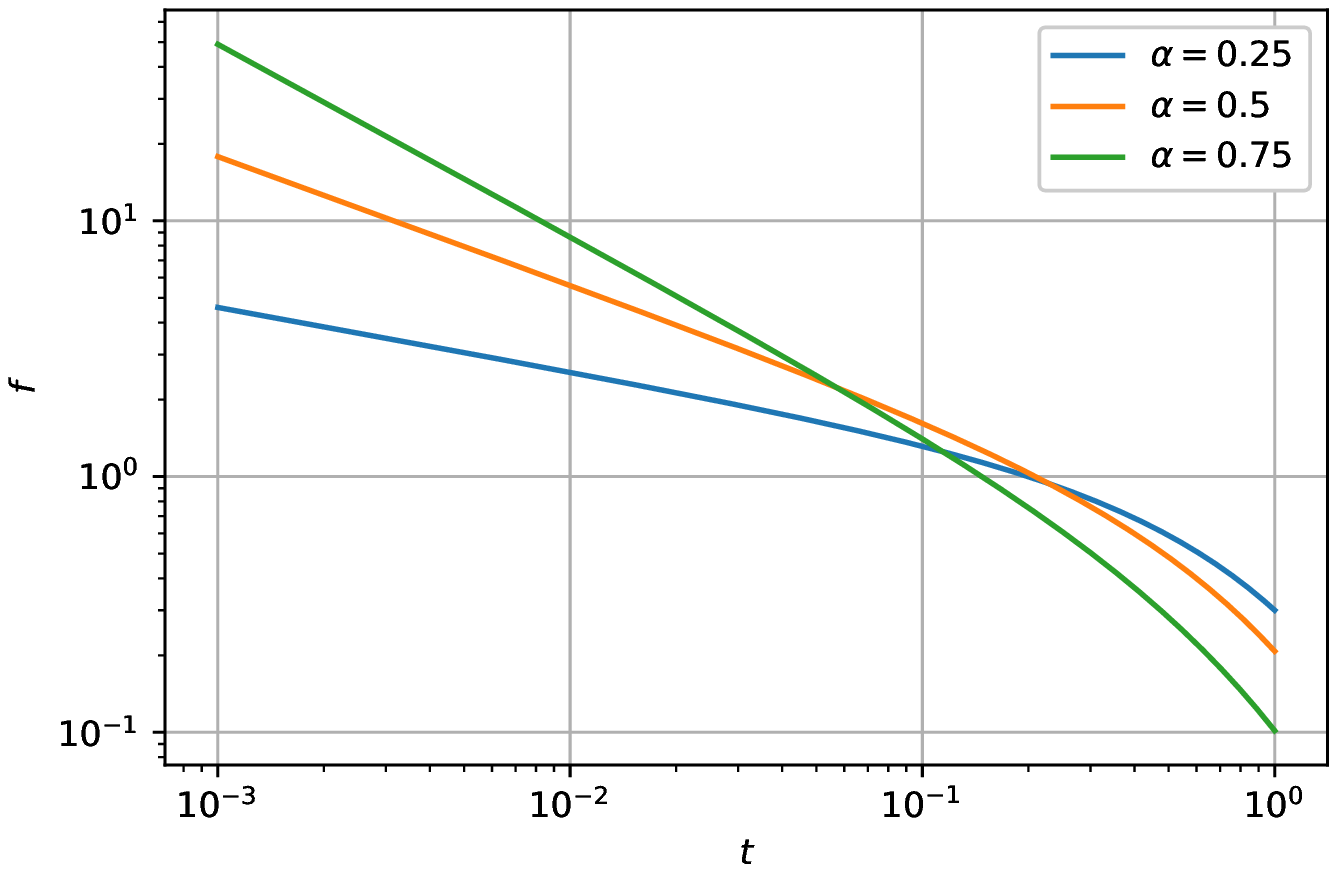} 
\caption{The kernel $k(t)$ at different values of $\delta$ ($\alpha = 0.25$, left) and at different values of $\alpha$ ($\delta = 1$, right).}
\label{f-1}
\end{figure}

\begin{figure}
\centering
\includegraphics[width=0.49\linewidth]{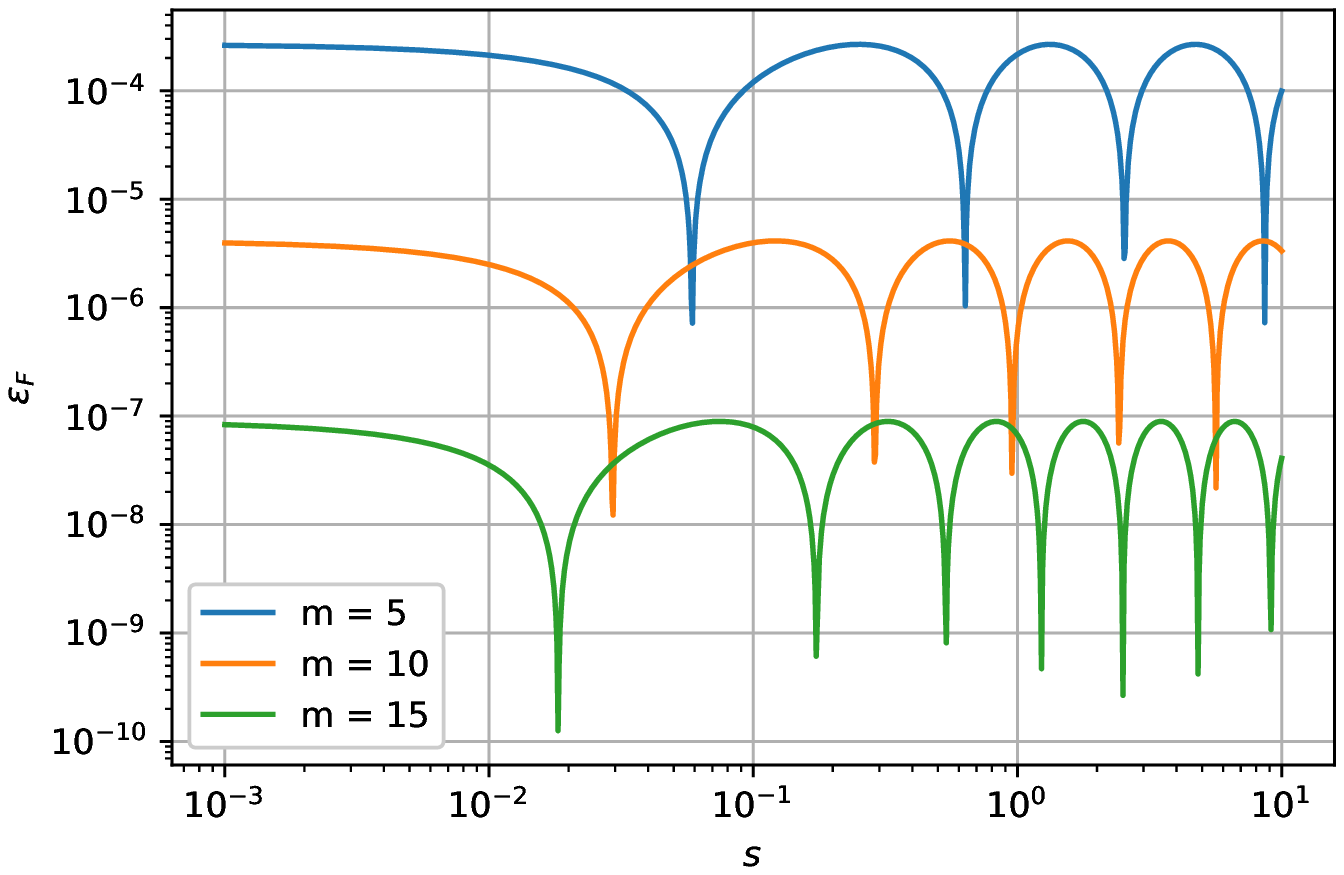} 
\includegraphics[width=0.49\linewidth]{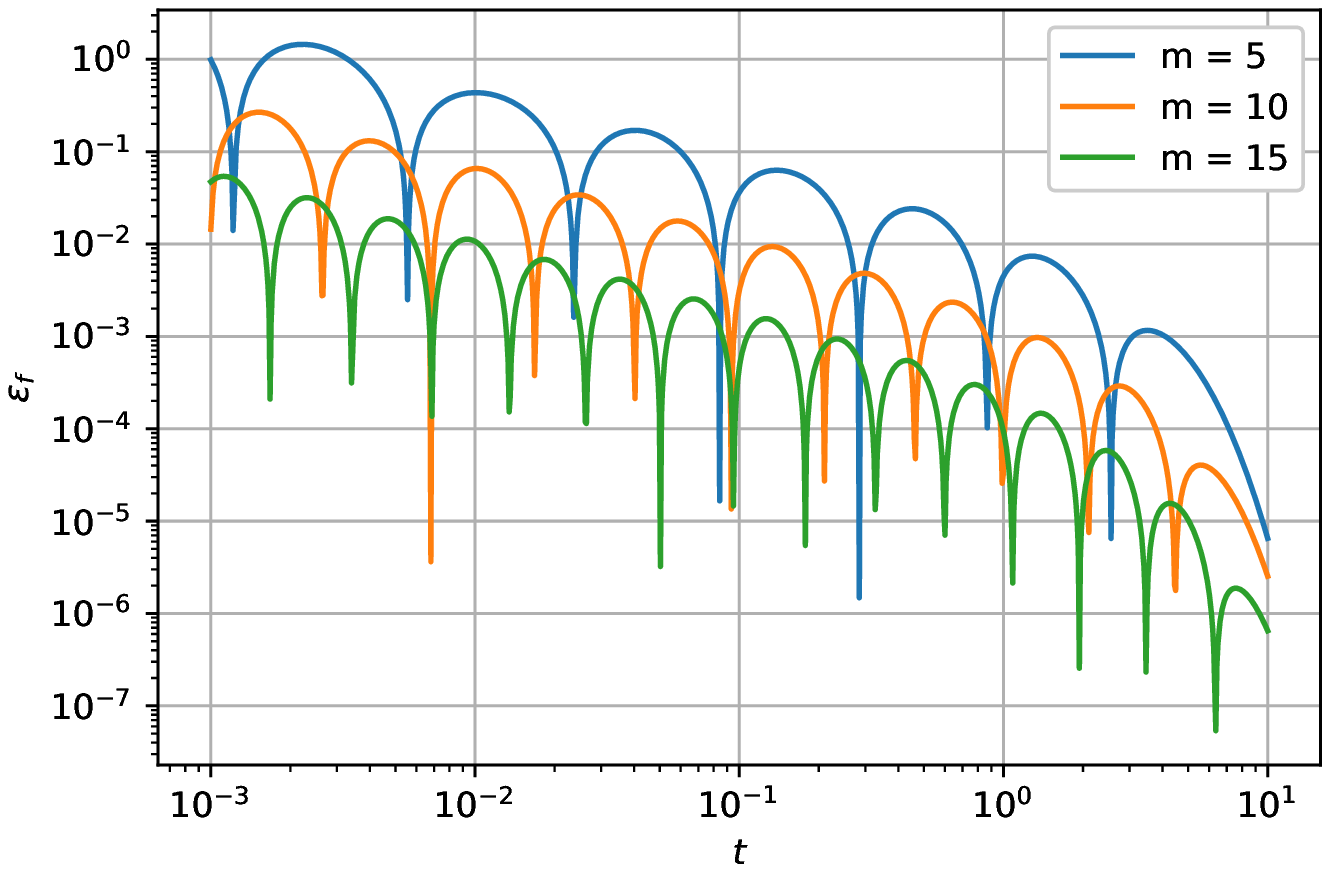} 
\caption{Approximation error $K(s)$ (left) and approximation error $k(t)$ (right) for different values of $m$ ($\alpha = 0.5, \delta=1$).}
\label{f-2}
\end{figure}

\begin{figure}
\centering
\includegraphics[width=0.49\linewidth]{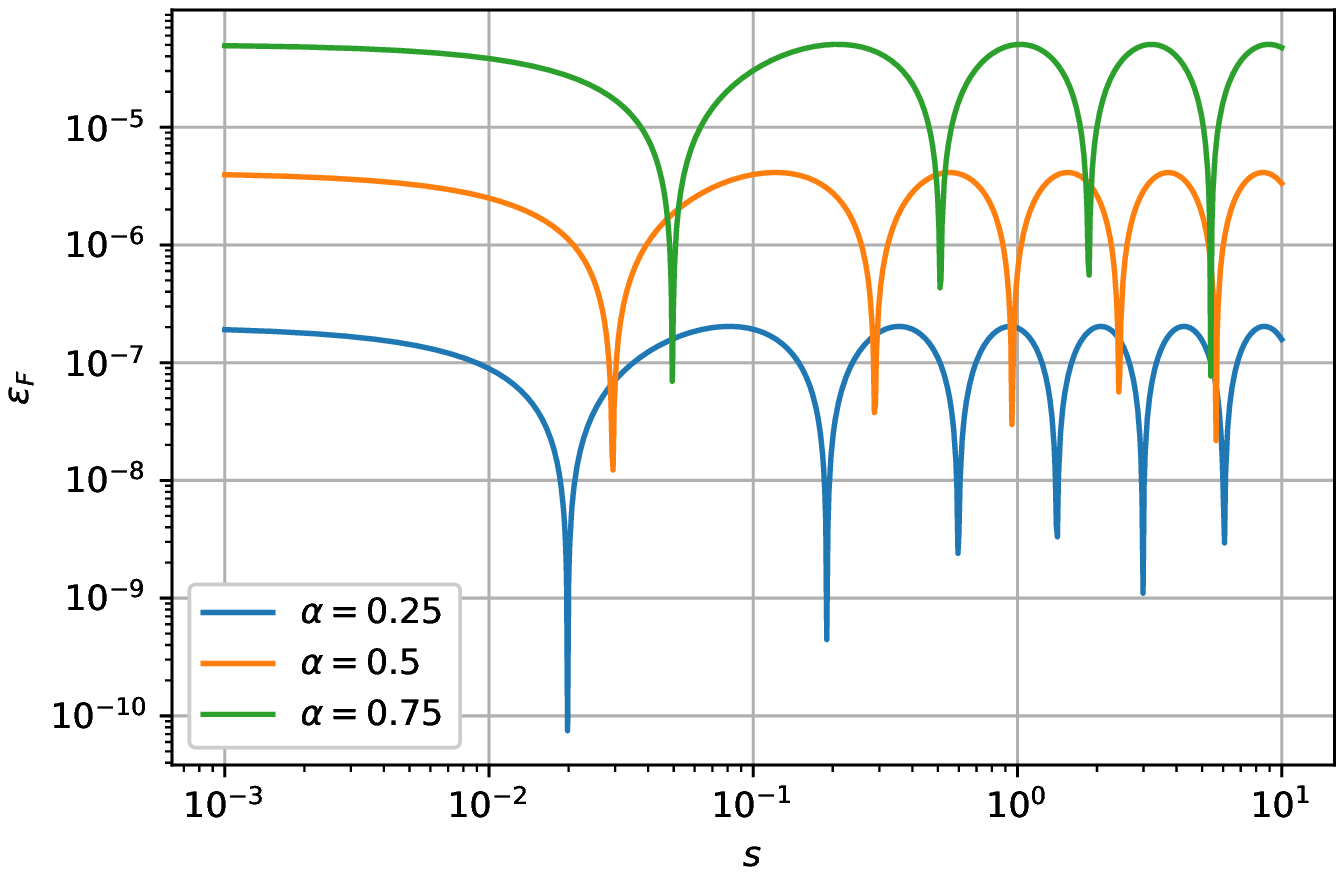} 
\includegraphics[width=0.49\linewidth]{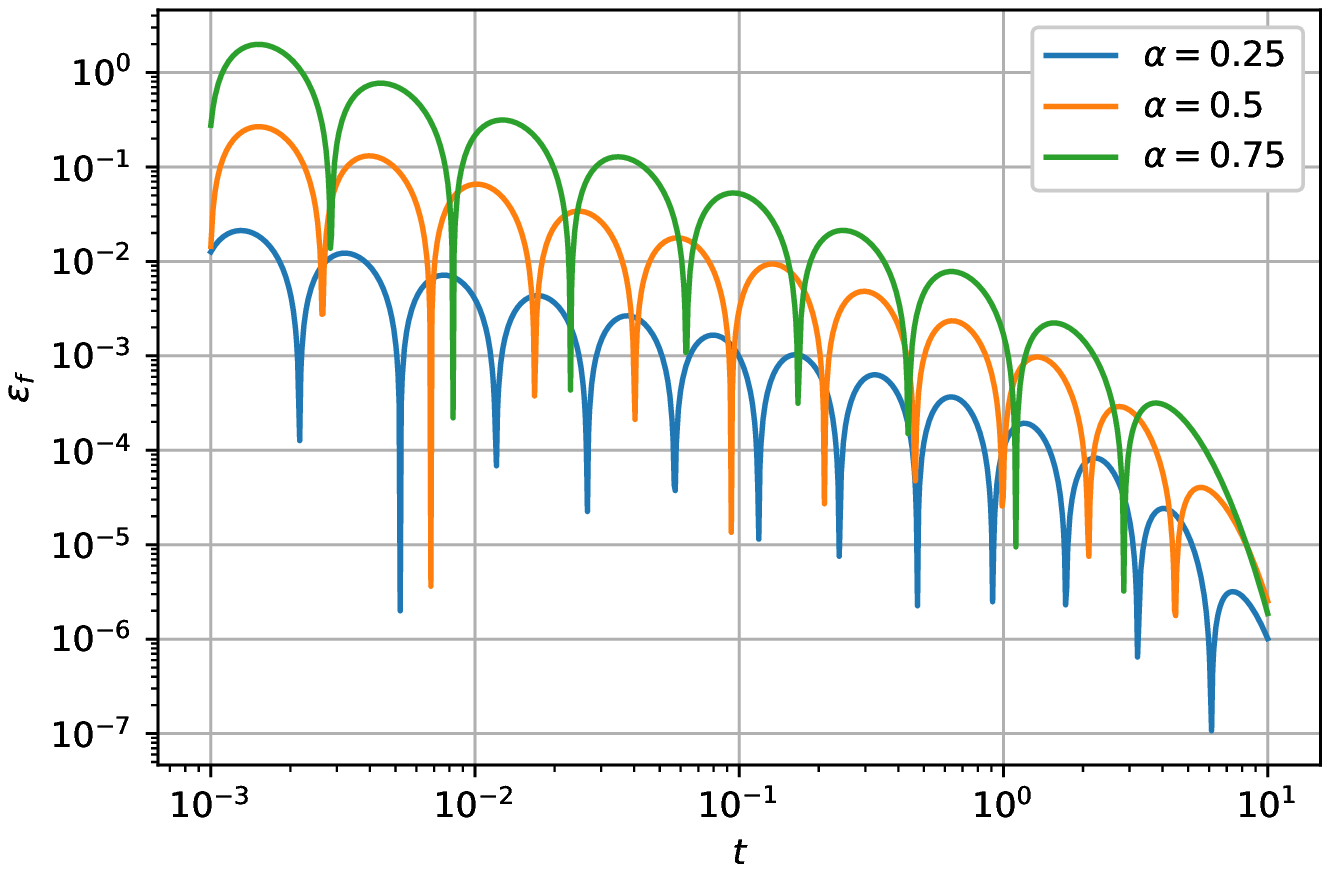} 
\caption{Approximation error $K(s)$ (left) and approximation error $k(t)$ (right) for different values of $\alpha$ ($m = 10, \delta=1$).}
\label{f-3}
\end{figure}

The kernel $k(t)$ for various values of the parameters $\alpha, \delta$ is shown in Fig.\ref{f-1}. 
As a base variant we consider the case $\alpha = 0.5, \delta=1$.
We observe a significant influence of $\alpha$, for small $t$ the influence of $\delta$ is insignificant.
A rational approximation was performed at $0 \leq s \leq 10^{3}$ for various $m$.
The approximation error $K(s)$ and $k(t)$ was estimated as follows:
\[
 \varepsilon_F(s) = |\accentset{\sim }{K}(s) - K(s)|,
 \quad \varepsilon_f(t) = |\accentset{\sim }{k}(t) - k(t)| .
\] 
The accuracy of approximations is illustrated by Fig.\ref{f-2}.  
The singularity of the kernel $k(t)$ leads to a significant drop in the approximation accuracy for small $t$.
The dependence of the approximation accuracy of $K(s)$ and $k(t)$ on the key parameter $\alpha$ at $m=10$ is shown in Fig.\ref{f-3}.  
We have used such approximations in approximate solutions to problems with memory of the time derivative of the solution.
The corresponding data on the coefficients $a_i, b_i, \ i =1,2,\ldots, m,$ are given in Table ~1. 

\begin{center}
\begin{table}[htp]
\label{tabl-1}
\caption{Parameters of approximation with $m=10$ for $k(t)$.}
\centering
\begin{tabular}{|ccccccc|}
\hline
 $\alpha $       &\multicolumn{2}{c}{0.25}    &\multicolumn{2}{c}{0.5}  &\multicolumn{2}{c|}{0.75}  \\
\hline
 $\gamma_2$       &\multicolumn{2}{c}{2.652102e-04 }    &\multicolumn{2}{c}{4.969023e-03}  &\multicolumn{2}{c|}{7.245547e-02}  \\
\hline
  $i$  & $a_i$  & $b_i$  &  $a_i$   &  $b_i$   &  $a_i$   &  $b_i$  \\

       1         &  5.521381e-01     &  1.020117e+00     &    2.819331e-01    &  1.047498e+00    &  1.104072e-01    &  1.083461e+00  \\
       2         &  2.880242e-01     &  1.366452e+00     &    3.375860e-01    &  1.485379e+00    &  2.289235e-01    &  1.631796e+00  \\
       3         &  2.752413e-01     &  2.383162e+00     &    4.623698e-01    &  2.727644e+00    &  4.363529e-01    &  3.155924e+00  \\
       4         &  3.105713e-01     &  4.913933e+00     &    6.873945e-01    &  5.845661e+00    &  8.492928e-01    &  7.015601e+00  \\
       5         &  3.787098e-01     &  1.118313e+01     &    1.072155e+00    &  1.364996e+01    &  1.692847e+00    &  1.675636e+01  \\
       6         &  4.824220e-01     &  2.712735e+01     &    1.730473e+00    &  3.361376e+01    &  3.462139e+00    &  4.174784e+01  \\
       7         &  6.372698e-01     &  6.947911e+01     &    2.904397e+00    &  8.670029e+01    &  7.385665e+00    &  1.081333e+02  \\
       8         &  8.891043e-01     &  1.907182e+02     &    5.242115e+00    &  2.388307e+02    &  1.729807e+01    &  2.986573e+02  \\
       9         &  1.408493e+00     &  5.988445e+02     &    1.131457e+01    &  7.591537e+02    &  5.160519e+01    &  9.649405e+02  \\
      10         &  3.276183e+00     &  2.794376e+03     &    4.224693e+01    &  3.797078e+03    &  3.306752e+02    &  5.301624e+03  \\
\hline
\end{tabular}
\end{table}
\end{center}

\begin{figure}
\centering
\includegraphics[width=0.49\linewidth]{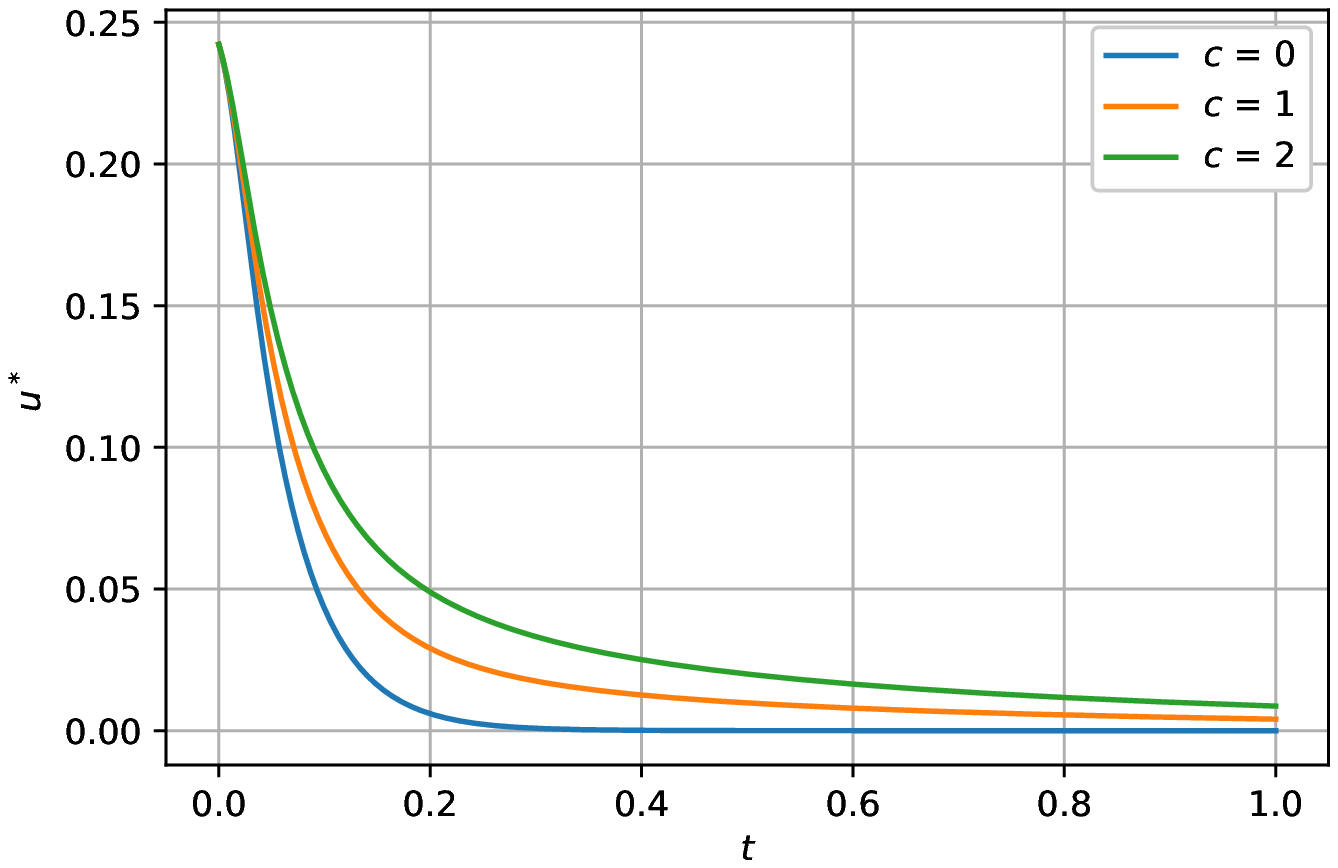} 
\includegraphics[width=0.49\linewidth]{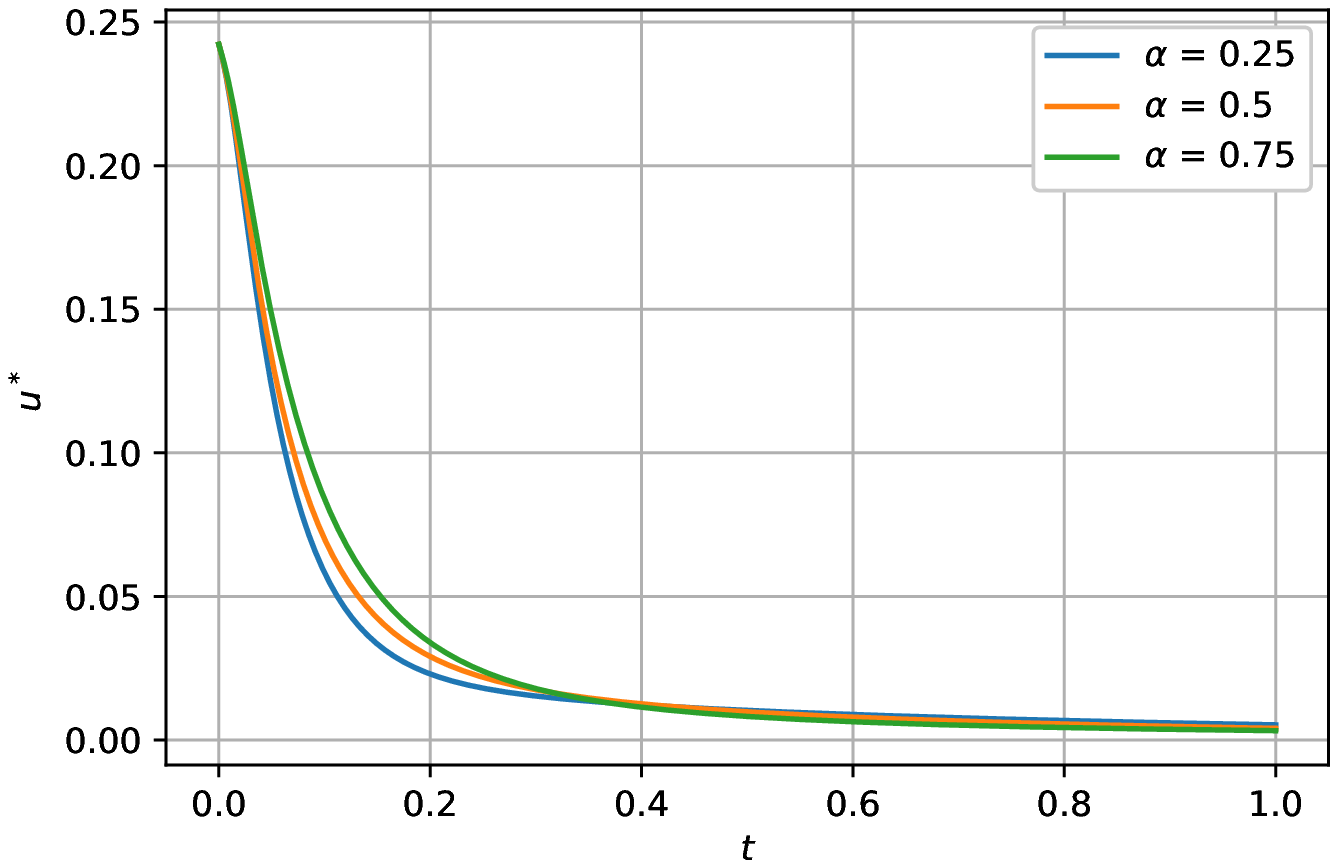} 
\caption{Solution $u^* = u(\bm x^*, t)$ at the point $\bm x^* = (0.5,0.5)$ for different values of the parameter $c, \alpha$: 
 $\alpha = 0.5$ (left) and $c = 1$ (right).} 
\label{f-4}
\end{figure}

\begin{figure}
\centering
\begin{minipage}{0.49\linewidth}
\centering
\includegraphics[width=\linewidth]{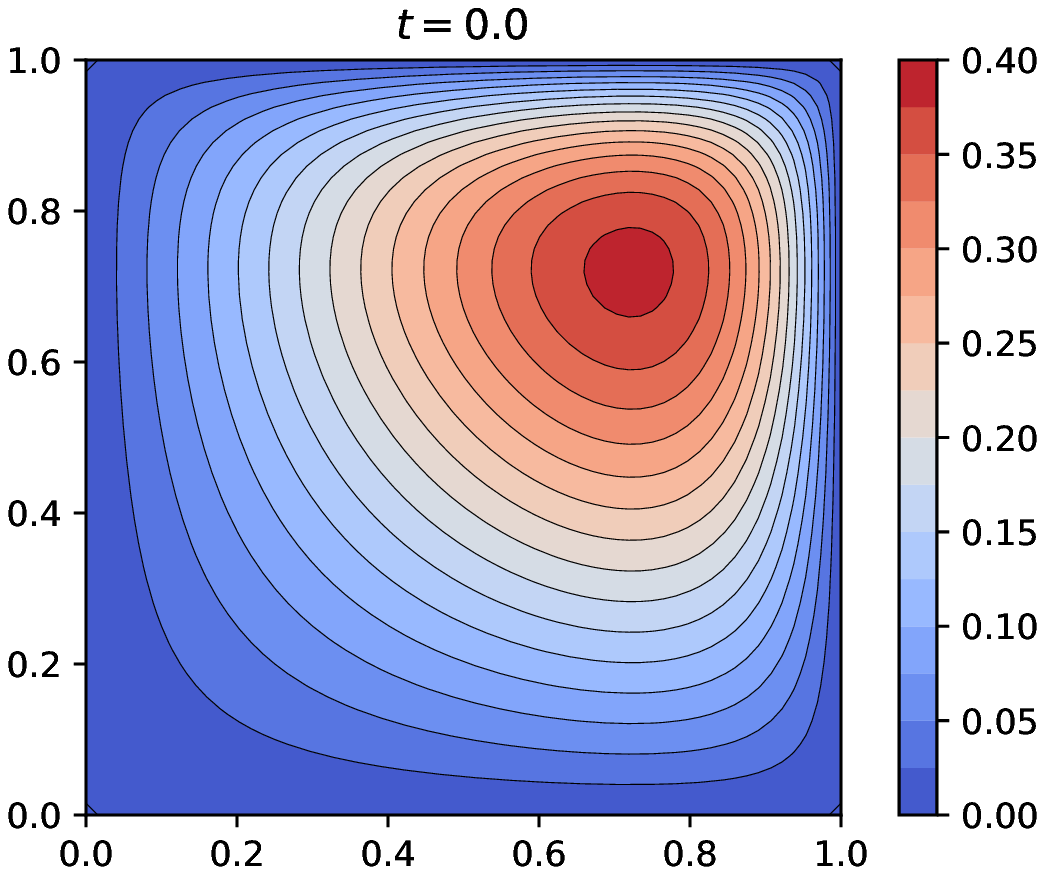} \\
\includegraphics[width=\linewidth]{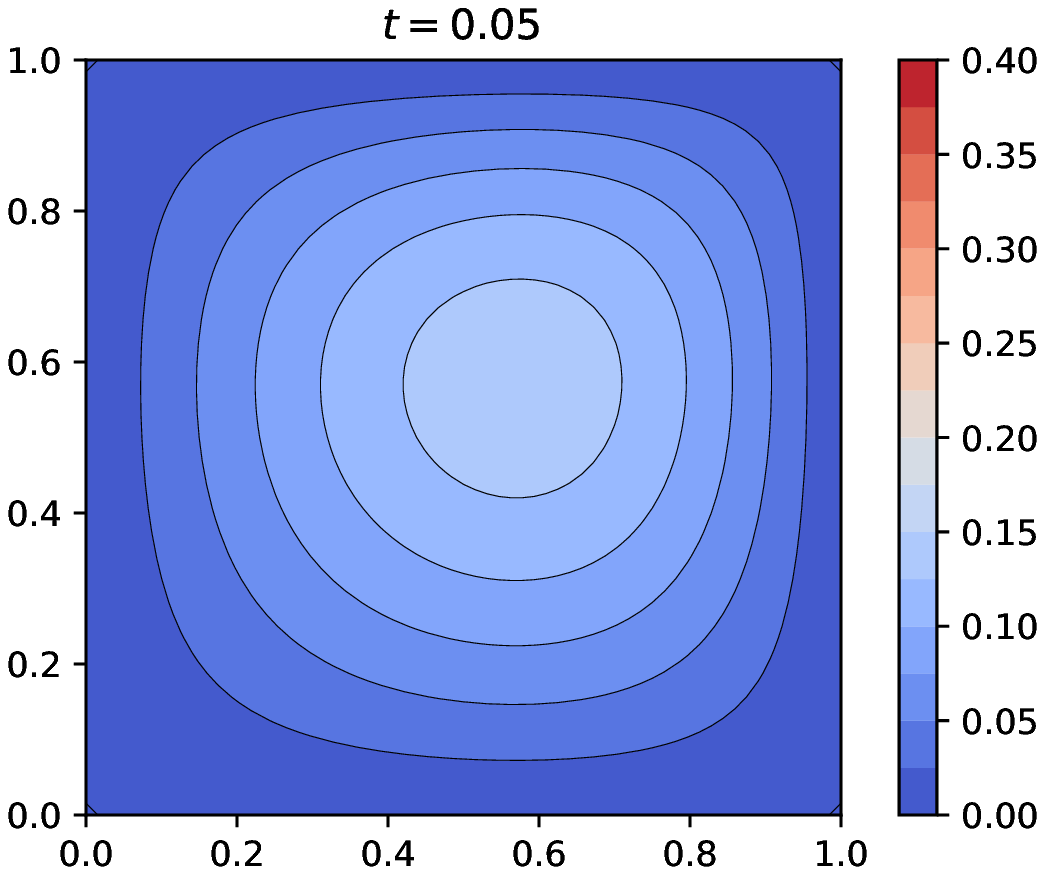} \\
\end{minipage}
\begin{minipage}{0.49\linewidth}
\centering
\includegraphics[width=\linewidth]{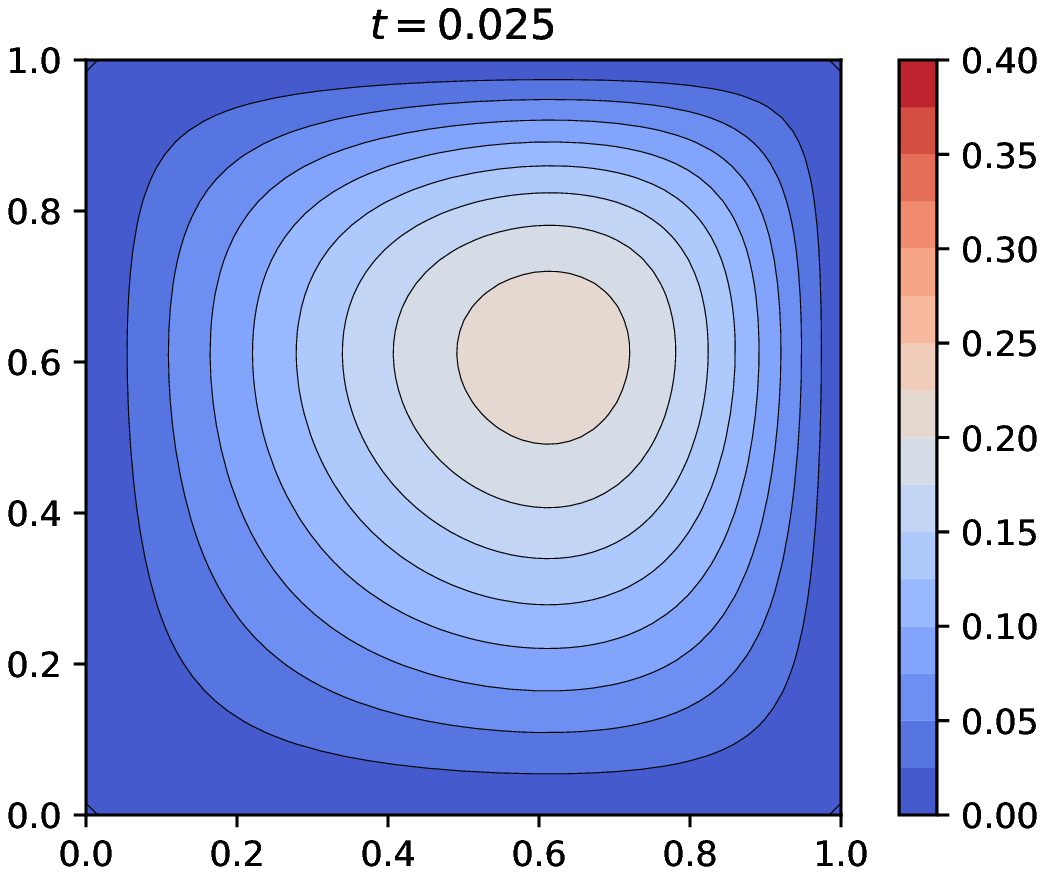} \\
\includegraphics[width=\linewidth]{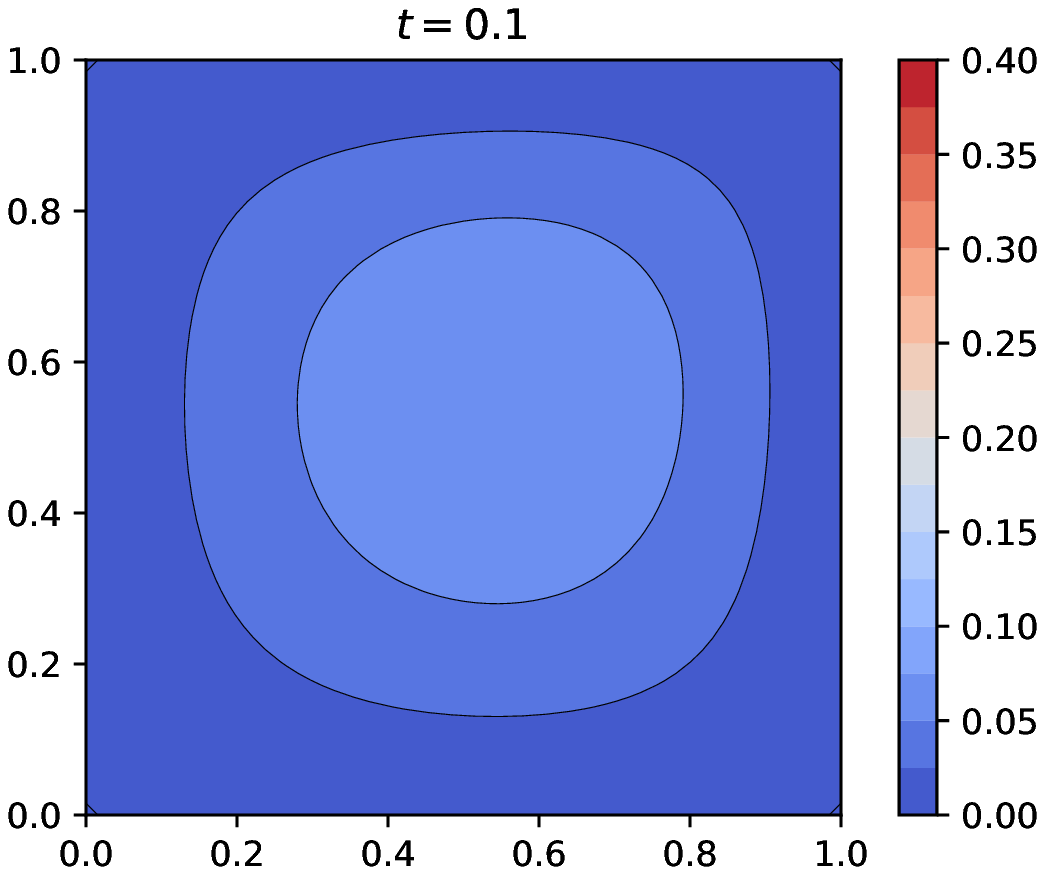} \\
\end{minipage}
\caption{Solution of the problem at separate points in time for $c=1, \alpha=0.5$.}
\label{f-5}
\end{figure}

We consider the problem of relaxation of the initial state of the system under study when
\[
 f(\bm x,t) = 0,
 \quad u_0(\bm x) = x_1 (1-x_1^6) x_2 (1-x_2^6) .
\] 
The calculations were performed on the spatial grid $N_1 = N_2 = 64$. 
The approximated solution at various time steps was compared with the solution on the detailed grid with $N=1000$, which was obtained using a second-order approximation scheme ($\sigma = 0.5$ in (\ref{4.5})--(\ref{4.7})). 
Fig.\ref{f-4} shows such a reference solution in the center of the computational domain for various values of the parameters $c, \alpha$. 
The presence of time derivative memory ($c > 0$) leads to a pronounced slowing of the decay rate of the solution.
The solution at separate moments of time for the test problem at $c=1, \alpha=0.5$ is shown in Fig.\ref{f-5}. 

\begin{figure}
\centering
\includegraphics[width=0.45\linewidth]{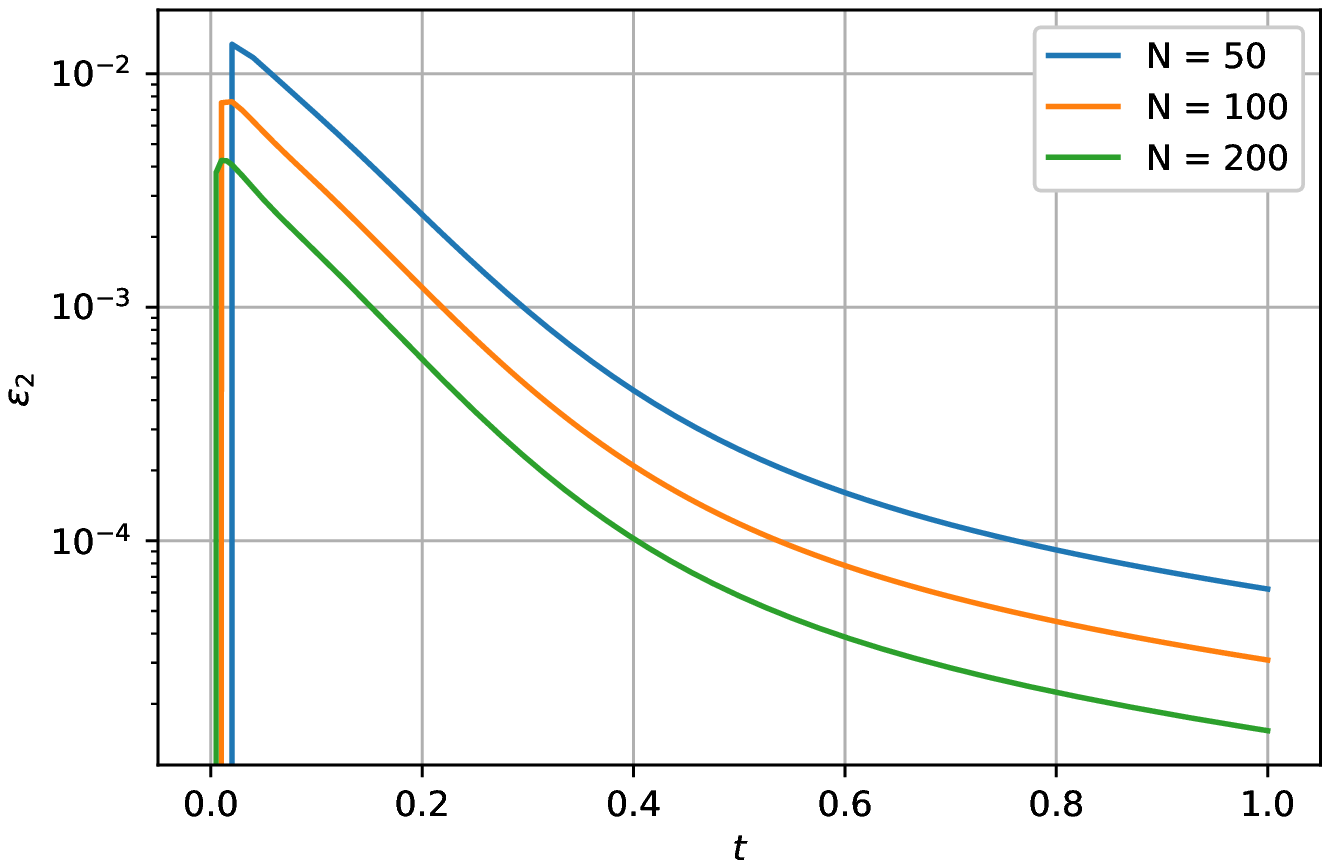} 
\includegraphics[width=0.45\linewidth]{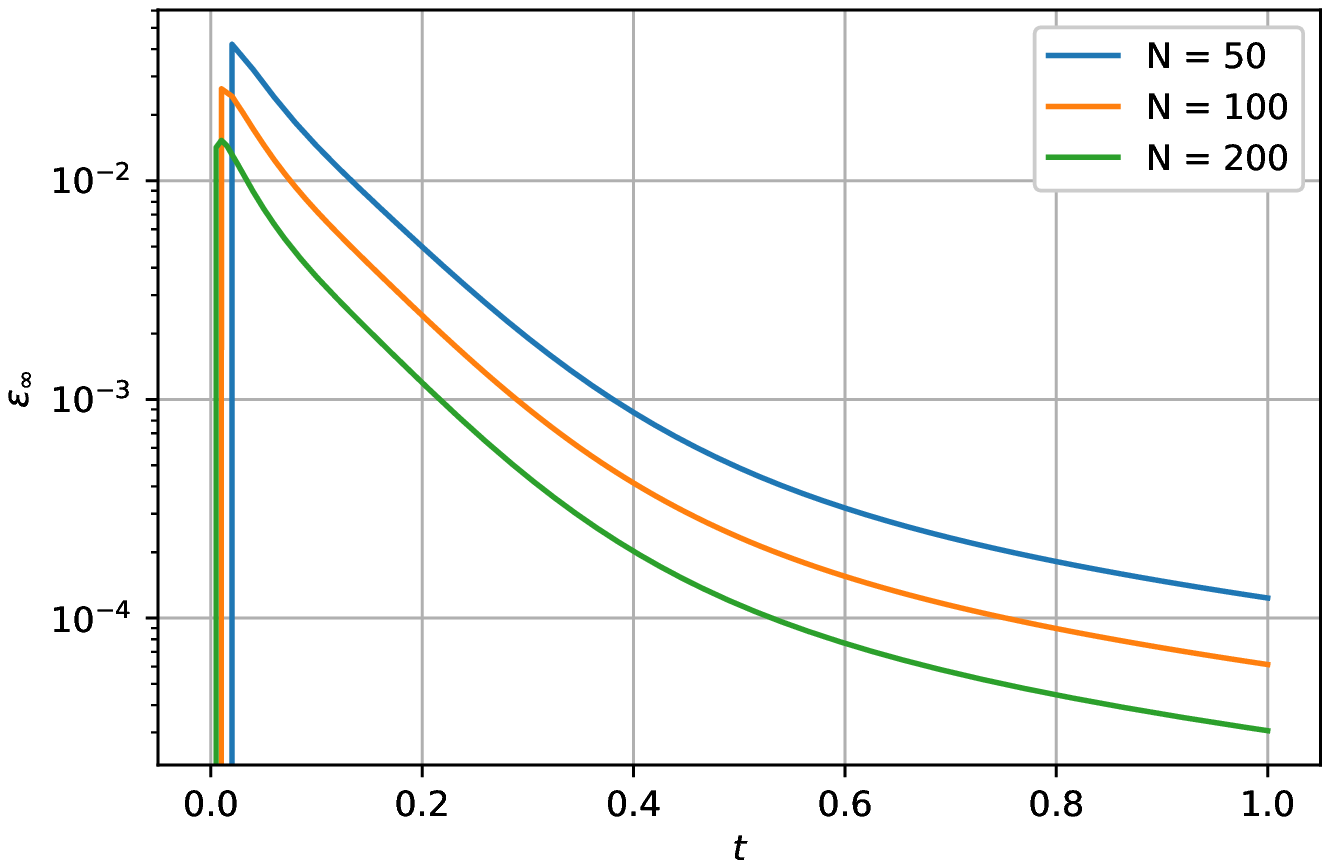}
\caption{Accuracy of the implicit Euler scheme ($\sigma = 1$) for the problem with $c=1, \alpha=0.5$.}
\label{f-6}
\end{figure}

\begin{figure}
\centering
\includegraphics[width=0.45\linewidth]{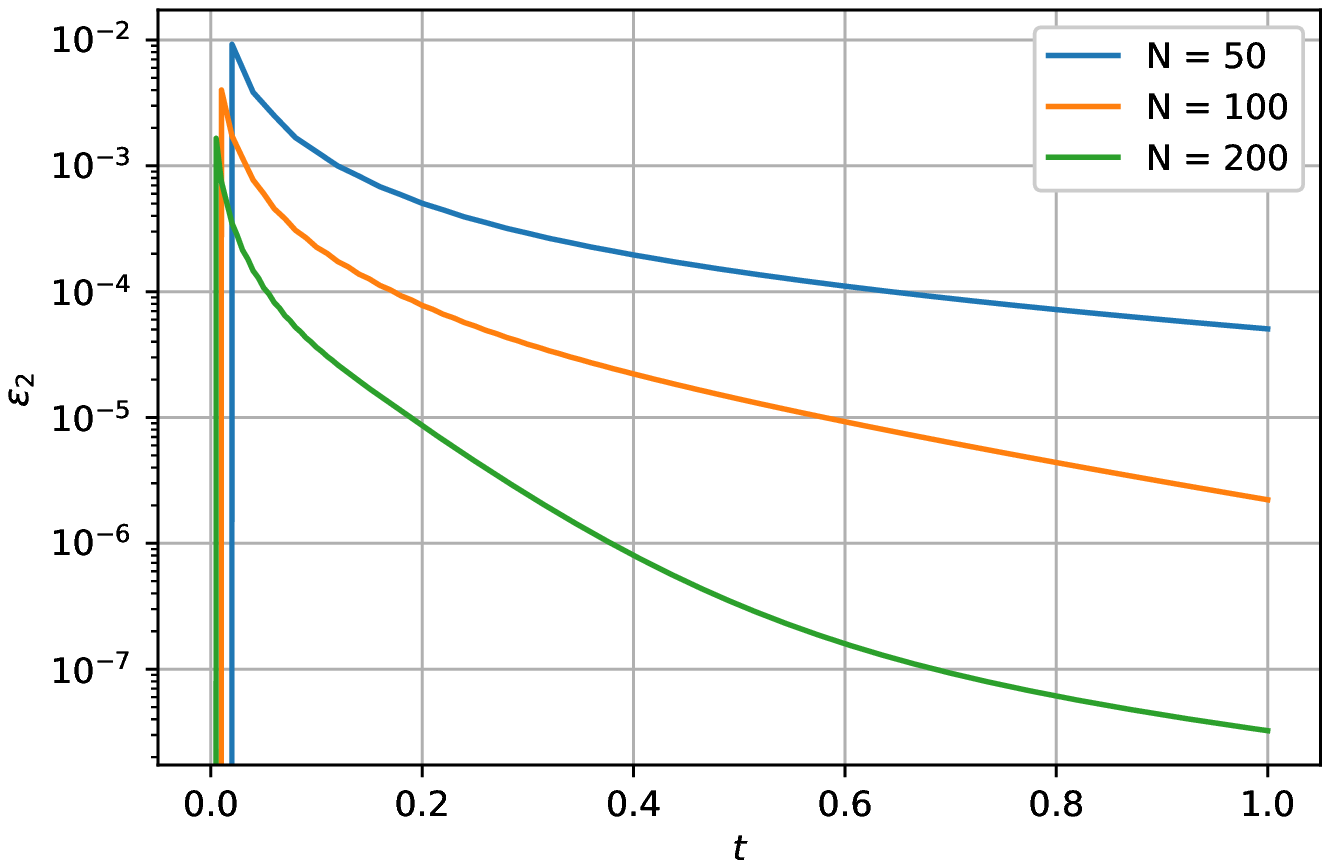} 
\includegraphics[width=0.45\linewidth]{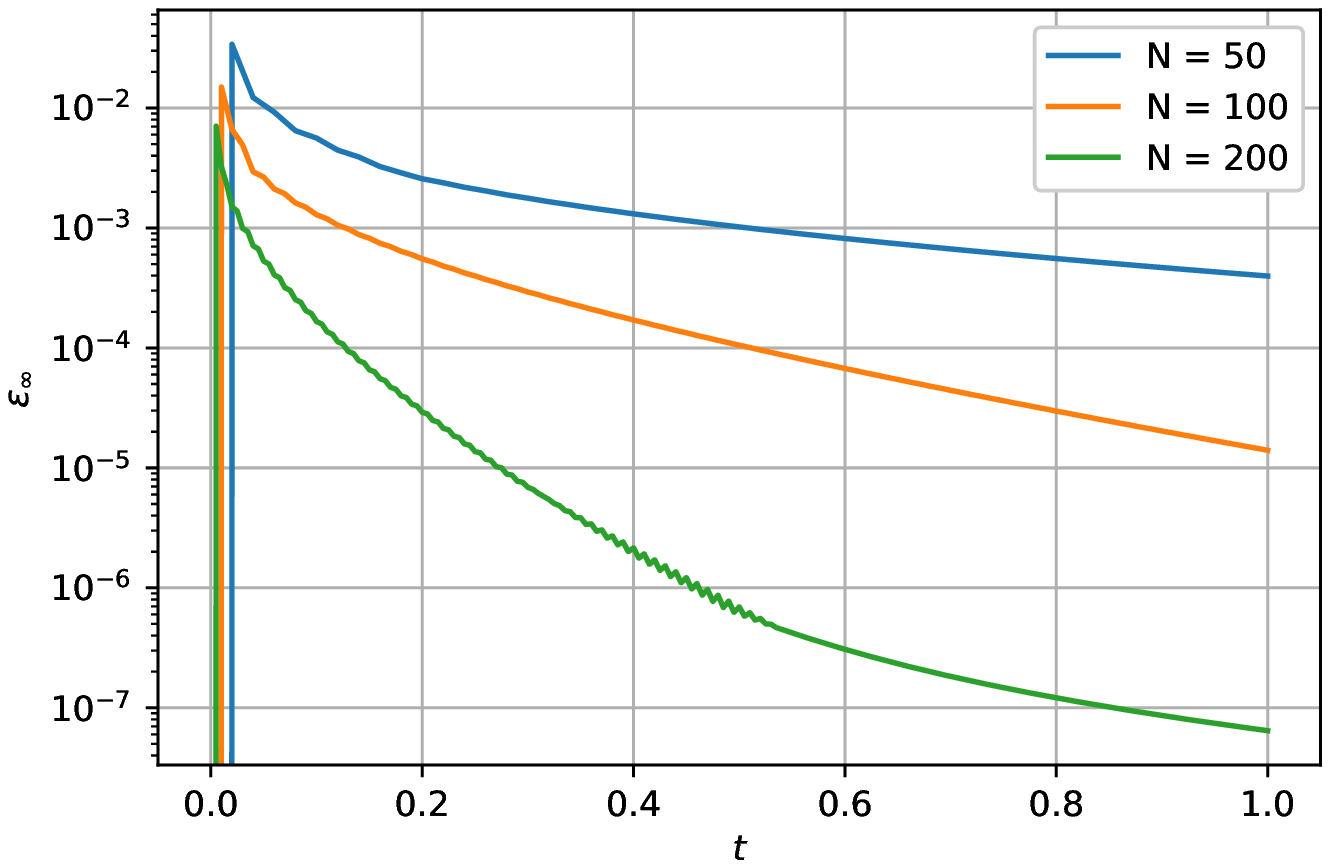} 
\caption{Accuracy of the symmetric scheme ($\sigma = 0.5$) for the problem with $c=1, \alpha=0.5$.}
\label{f-7}
\end{figure}

The accuracy of the solution of the problem with memory is estimated by the absolute discrepancy at individual points in time:
\[
 \varepsilon_2(t^n) = \|y (\bm x,t^n) - \bar{y}^n (\bm x,t^n)\| ,
 \quad \varepsilon_\infty (t^n) = \max_{\bm x \in \omega} |y (\bm x,t^n) - \bar{y}^n(\bm x)| ,
 \quad n = 0, \ldots, N ,
\] 
where $\bar{y}$ is the reference solution. 
The accuracy when using the implicit Euler scheme is shown in Fig.\ref{f-6}.
Similar data for the symmetric scheme are shown in Fig.\ref{f-7}. 
The calculated data are consistent with the above theoretical considerations about the accuracy of two-level schemes with weights (\ref{4.5})--(\ref{4.7}) for not very small $t$. 
At the initial stage, the solution changes most strongly, which is reflected in a drop in the accuracy of the approximate solution at $t \ll 1$. 

\section{Conclusions}\label{sec:6} 

\begin{enumerate}
\item We considered the Cauchy problem for a first-order integrodifferential equation with the memory of the time derivative of the solution.
The kernel is assumed to be difference, and the operators of the equation are self-adjoint and positive definite in a finite-dimensional Hilbert space.
We have established that the solution is stable concerning the right-hand side and initial conditions under the usual assumption that the kernel is positive definite.
The problems of numerical solution of such problems are mainly related to the necessity to operate with the solution for all previous moments.
 \item A well-known approach based on the approximation of a difference kernel by a sum of exponentials is used for an approximate solution of the formulated evolutionary problem with memory.
In this case, we pass from a nonlocal problem in time to a local problem for a system of weakly coupled evolution equations with additional ordinary differential equations for auxiliary functions.
We give a priori estimates for the solution of the Cauchy problem for this system of evolutionary equations.
 \item The numerical solution uses standard two-level time approximations.
The unconditional stability of two-level schemes with weights under standard constraints on weights is proved.
The transition to a new level in time is provided by solving the usual problem for the approximate solution itself and explicitly recomputing the auxiliary functions.
 \item We complemented the theoretical consideration with the data of numerical solution of the model two-dimensional problem.
Numerical approximation of the difference kernel by the sum of exponents is performed based on rational approximation, which is applied to the Laplace transform for the kernel. The influence of various parameters on the approximate solution of the problem is investigated.
\end{enumerate} 


\end{document}